\documentclass[10pt, a4paper]{article}
\usepackage{amscd}
\usepackage{amsfonts}
\usepackage{caption}
\begin{document}

\newtheorem{theorem}{Theorem}
\newtheorem{corollary}[theorem]{Corollary}
\newtheorem{definition}[theorem]{Definition}
\newtheorem{lemma}[theorem]{Lemma}
\newtheorem{proposition}[theorem]{Proposition}
\newtheorem{remark}[theorem]{Remark}
\newtheorem{example}[theorem]{Example}
\newtheorem{notation}[theorem]{Notation}
\def\Qed{\hfill\raisebox{.6ex}{\framebox[2.5mm]{}}\\[.15in]}

\title{Involutions on surfaces with\\ $p_g=q=0$ and $K^2=3$}

\author{Carlos Rito}

\date{}
\pagestyle{myheadings}
\maketitle
\setcounter{page}{1}

\begin{abstract}
We study surfaces of general type $S$ with $p_g=0$ and $K^2=3$ having an involution $i$ such that the bicanonical map of $S$ is not composed with $i$. It is shown that, if $S/i$ is not rational, then $S/i$ is birational to an Enriques surface or it has Kodaira dimension $1$ and the possibilities for the ramification divisor of the covering map $S\rightarrow S/i$ are described. We also show that these two cases do occur, providing an example. In this example $S$ has a hyperelliptic fibration of genus $3$ and the bicanonical map of $S$ is of degree $2$ onto a rational surface.

\end{abstract}

\section{Introduction}

Minimal surfaces $S$ of general type with $p_g=q=0$ have been studied by several authors in the last years, but a classification is still missing. For these surfaces the canonical divisor $K$ satisfies $1\leq K^2\leq 9$ and there are examples for all values of $K^2$ (see {\it e.g.} \cite{BHPV}).
The study of the bicanonical map $\phi_2$ of $S,$ and in particular the case where $\phi_2$ is composed with an involution of $S,$ has also provided some examples ({\it cf.} \cite{CFM}, \cite{MP1}, \cite{MP2}, \cite{MP4}, \cite{MP5}, \cite{MP6}).

For the case $K^2=3,$ there are examples with bicanonical map of degree $2$ onto a nodal Enriques surface (see \cite{MP3}, \cite{MP5}) and with bicanonical map of degree $4$ onto a rational surface (see \cite{Bu}, \cite{In}, \cite{Ca}, \cite{Ke}, \cite{Na}). In these constructions with $\deg(\phi_2)=4$ the surface $S$ has involutions $i_j,$ $j=1,2,3,$ such that $\phi_2$ is composed with $i_j$ and such that $S/{i_j}$ is birational to an Enriques surface or $S/{i_j}$ is a rational surface.

There are also the constructions given in \cite{PPS1} and \cite{PPS2}, obtained using $\mathbb Q$-Gorenstein smoothing theory,
and the recent construction given in \cite{BP},
but we have no information about the bicanonical map or the existence of an involution in these cases.

In this paper we want to study the case where $K_S^2=3$ and $S$ has an involution $i$ such that the bicanonical map of $S$ is not composed with $i.$ We show that, if $S/i$ is not rational, then $S/i$ is birational to an Enriques surface or it has Kodaira dimension $1$ and we describe the possibilities for the ramification divisor of the covering map $S\rightarrow S/i.$ We also show that these two cases do occur, providing an example. In this example $S$ has a hyperelliptic genus $3$ fibration and the bicanonical map of $S$ is of degree $2$ onto a rational surface.

The paper is organized as follows. In Sections \ref{GenFacts} and \ref{NumRestr} we recall some facts about involutions on surfaces and about the possibilities for the {\em branch locus} (the projection of the ramification divisor) in the quotient surface $S/i.$ This is used to prove our main results in Section \ref{Poss}, Theorems \ref{main} and \ref{components}. Section \ref{example} contains the construction of an example, which is obtained as a $\mathbb Z_2^2$ cover of $\mathbb P^2.$ The ramification divisor of this covering is computed using the {\em Computational Algebra System Magma} (\cite{BCP}). The corresponding code lines are given in the Appendix.

\bigskip
\noindent{\bf Notation}

We work over the complex numbers; all varieties are assumed to be projective algebraic.

An {\em involution} of a surface $S$ is an
automorphism of $S$ of order 2. We say that a map is {\em composed with
an involution} $i$ of $S$ if it factors through the double cover $S\rightarrow
S/i.$

An {\em $(-2)$-curve} $N$ on a surface is a curve isomorphic to $\mathbb P^1$ such that $N^2=-2$.

An $(m_1,m_2,\ldots)$-point of a curve, or point of type $(m_1,m_2,\ldots),$ is a singular point of multiplicity $m_1,$ which resolves to a point of multiplicity $m_2$ after one blow-up, etc.

The rest of the notation is standard in Algebraic Geometry.\\

\bigskip
\noindent{\bf Acknowledgements}

The author wishes to thank Margarida Mendes Lopes for all the support.
He is a member of the Mathematics Center of the Universidade de Tr\'as-os-Montes e Alto Douro and is a collaborator of the Center for Mathematical Analysis, Geometry and Dynamical Systems of Instituto Superior T\'ecnico, Universidade T\' ecnica de Lisboa.
This research was partially supported by FCT (Portugal) through Project PTDC/MAT/099275/2008.

\section{General facts on involutions}\label{GenFacts}

The following is according to \cite{CM}.\\
Let $S$ be a smooth minimal surface of general type with an
involution $i.$ Since $S$ is minimal of general type, this
involution is biregular. The fixed locus of $i$ is the union of a
smooth curve $R''$ (possibly empty) and of $t\geq 0$ isolated points
$P_1,\ldots,P_t.$ Let $S/i$ be the quotient of $S$ by $i$ and
$p:S\rightarrow S/i$ be the projection onto the quotient. The
surface $S/i$ has nodes at the points $Q_i:=p(P_i),$ $i=1,\ldots,t,$
and is smooth elsewhere. If $R''\not=\emptyset,$ the image via $p$
of $R''$ is a smooth curve $B''$ not containing the singular points
$Q_i,$ $i=1,\ldots,t.$ Let now $h:V\rightarrow S$ be the blow-up of
$S$ at $P_1,\ldots,P_t$ and set $R'=h^*(R'').$ The involution $i$
induces a biregular involution $\widetilde{i}$ on $V$ whose fixed
locus is $R:=R'+\sum_1^t h^{-1}(P_i).$ The quotient
$W:=V/\widetilde{i}$ is smooth and one has a commutative diagram:
$$
\begin{CD}\ V@>h>>S\\ @V\pi VV  @VV p V\\ W@>g >> S/i
\end{CD}
$$
where $\pi:V\rightarrow W$ is the projection onto the quotient and
$g:W\rightarrow S/i$ is the minimal desingularization map. Notice
that $$A_i:=g^{-1}(Q_i),\ \ i=1,\ldots,t,$$ are $(-2)$-curves and
$\pi^*(A_i)=2\cdot h^{-1}(P_i).$ 

Set $B':=g^*(B'').$ Since $\pi$ is
a double cover with branch locus $B'+\sum_1^t A_i,$ it is determined
by a line bundle $L$ on $W$ such that $$2L\equiv B:=B'+\sum_1^t A_i.$$
It is well known that ({\it cf.} \cite[V. 22]{BHPV}):
$$p_g(S)=p_g(V)=p_g(W)+h^0(W,\mathcal{O}_W(K_W+L)),$$
$$q(S)=q(V)=q(W)+h^1(W,\mathcal{O}_W(K_W+L)),$$
$$K_S^2-t=K_V^2=2(K_W+L)^2$$
and
$$\chi(\mathcal{O}_S)=\chi(\mathcal{O}_V)=2\chi(\mathcal{O}_W)+\frac{1}{2}L(K_W+L).$$
\begin{lemma}[\cite{CM}, \cite{CCM}]
The bicanonical map $\phi_2$ of $S$ (given by $|2K_S|$) is composed with $i$ if and only if $h^0(W,\mathcal{O}_W(2K_W+L))=0$.
\end{lemma}

\section{Numerical restrictions}\label{NumRestr}

Let $P$ be a minimal model of the resolution $W$ of $S/i$ and
$\rho:W\rightarrow P$ be the corresponding projection. Denote by
$\overline{B}$ the projection $\rho(B)$ and by $\delta$ the
"projection" of $L.$
\begin{remark}\label{CanRes}
If $\overline B$ is singular, there are exceptional divisors $E_i$ and numbers $r_i\in 2\mathbb N$ such that
$$
\begin{array}{l}
E_i^2=-1,\\
K_W\equiv\rho^*(K_P)+\sum E_i,\\
2L\equiv B=\rho^*(\overline{B})-\sum r_iE_i\equiv\rho^*(2\delta)-\sum r_iE_i.
\end{array}
$$
\end{remark}
The next result follows from Propositions 2, 3 a) and 4 b) of \cite{Ri1}.
\begin{proposition}[{\it cf.} \cite{CM}, \cite{Ri1}]\label{NumR}
Let $S$ be a smooth minimal surface of general type with $p_g=0$ and $K^2=3$ having an involution $i.$
With the previous notation, we have:
\begin{description}
  \item[a)] $K_P(K_P+\delta)+\frac{1}{2}\sum(r_i-2)=h^0(W,\mathcal{O}_W(2K_W+L));$
  \item[b)] $\delta^2=-2K_P^2-3K_P\delta+\frac{1}{4}\sum(r_i-2)(r_i-4)+2h^0(W,\mathcal{O}_W(2K_W+L))-2;$
  \item[c)] the number of isolated fixed points of $i$ is $t=7-2h^0(W,\mathcal{O}_W(2K_W+L));$
  \item[d)] $K_W^2\geq 2h^0(W,\mathcal{O}_W(2K_W+L))-4.$
\end{description}
\end{proposition}

\section{Possibilities}\label{Poss}

If $p_g(S)=0$ and the bicanonical map of $S$ is composed with the involution $i,$ then it is known that $S/i$ is birational to an Enriques surface or $S/i$ is a rational surface ({\it cf.} \cite{MP3}, \cite{MP5}). This follows easily from Proposition \ref{NumR}, a), b): we have $K_P(K_P+\delta)+\frac{1}{2}\sum(r_i-2)=0,$ thus $K_P$ nef implies $K_P^2=K_P\delta=0.$
Hence $P$ is birational to an Enriques surface or it has Kodaira dimension $1.$ In this last case $K_P\delta=0$ implies the existence of an elliptic fibration in $S,$ which is impossible because $S$ is of general type. 

Consider the branch divisor $B=B'+\sum_1^t A_i\subset W$ as above and let $\overline B,$ $\overline {B'}$ be the projection of $B,$ $B'$ on the minimal model $P$ of $W.$
We have the following:

\begin{theorem}\label{main}

Let $S$ be a smooth minimal surface of general type with $p_g=0$ and $K^2=3$ having an involution $i$ such that the bicanonical map of $S$ is not composed with $i.$

Then the number of isolated fixed points of $i$ is $t=5$ and, if $S/i$ is not rational, one of the following holds:
\begin{description}
  \item[a)] $P$ is an Enriques surface and
  	 \begin{description}
	   \item[(i)] $\overline {B'}^2=10,$ $\overline {B'}$ has a quadruple point and at most one double point $($thus $p_a(B')=0$ or $-1)$, or
	   \item[(ii)] $\overline {B'}^2=8,$ $\overline {B'}$ has a $(3,3)$-point and no other singularities $($thus $p_a(B')=-1)$;
	 \end{description}
  \item[b)] ${\rm Kod}(P)=1,$ $p_g(P)=q(P)=0,$ $\overline {B'}^2=-2,$ $p_a(\overline {B'})=1$ and $\overline {B'}$ has at most two double points.
\end{description}
Moreover, cases {\rm a) (i)} and {\rm b)} do occur; there is an example with bicanonical map of degree $2$ onto a rational surface.
\end{theorem}
{\bf Proof :}\\
Proposition 3 c) of \cite{Ri1} gives $h^0(W,\mathcal{O}_W(2K_W+L))\leq 1.$ Since $\phi_2$ is not composed with $i,$ we have $h^0(W,\mathcal{O}_W(2K_W+L))=1.$ Then, from Proposition \ref{NumR}, $t=5$ and $K_W^2\geq -2.$
\\
\\
{\bf Case 1:} ${\rm Kod}(P)=0.$\\
We have $p_g(P)\leq p_g(S)$ and $q(P)\leq q(S).$ Thus $p_g(P)=q(P)=0$ and then, from the classification of surfaces (see {\it e.g.} \cite{Be} or \cite{BHPV}), $P$ is an Enriques surface.
We obtain from Proposition \ref{NumR} that $\sum(r_i-2)=2$ and $\overline B^2=(2\delta)^2=0.$

Moreover, since $K_WA_i=0,$ each $(-2)$-curve $A_i\subset B$ is contracted to a singular point of $\overline{B'}$ or is mapped onto a $(-2)$-curve of the Enriques surface $P.$
\\
\\
{\bf Case 2:} ${\rm Kod}(P)=1.$\\
In this case $K_P$ is numerically equivalent to a rational multiple of a fibre of an elliptic fibration of $P$
(see {\it e.g.} \cite[V. 12]{BHPV}). 
This implies $K_P\overline B\ne 0,$ because otherwise $\overline B$ is contained in the elliptic fibration of $P$ and then $S$ has an elliptic fibration, which is impossible since $S$ is of general type. From Proposition \ref{NumR}, a) and b) we get
$\sum(r_i-2)=0,$ $K_P\overline B=2K_P\delta=2$ and $\overline B^2=(2\delta)^2=-12.$
\\
\\
{\bf Case 3:} ${\rm Kod}(P)=2.$
\\\\
{\bf Claim :} If $K_P\overline B=0,$ then $\overline B=B$ is a disjoint union of $(-2)$-curves.\\
{\em Proof}\ : Since $P$ is minimal of general type, $K_P$ is nef and big and then every component of $\overline B$ is a $(-2)$-curve and the intersection form on the components of the reduced effective divisor $\overline B$ is negative definite by the Algebraic Index Theorem (see {\rm e.g.} \cite[IV. 2.16]{BHPV}). The claim is true if each connected component $C$ of $\overline B$ is irreducible.
If $C$ is not irreducible, there is one component $\theta$ of $C$ such that $\theta(C-\theta)=1$ and this implies that $B$ has a $(-3)$-curve, contradicting $B\equiv 0\ ({\rm mod\ 2})$.\ $\diamondsuit$
\\

Since $P$ is of general type, $K_P^2\geq 1.$ Hence Proposition \ref{NumR} implies $K_P\delta=0,$  $\sum(r_i-2)=0$ and $\delta^2=-2.$ Therefore $\overline B^2=-8$ and then $B$ is a disjoint union of four $(-2)$-curves. But Proposition \ref{NumR} c) gives $t=5\ne 4$.
\\

An example for a) (i) and b) is given in Section \ref{example}. \Qed

In the conditions of Theorem \ref{main}, $S/i$ is a rational surface, or $S/i$ is birational to an Enriques surface or ${\rm Kod}(S/i)=1.$
We have no example for $S/i$ rational (and $\phi_2$ not composed with $i$) but there is at least one possibility that may occur: $S$ is the smooth minimal model of a double cover of $\mathbb P^2$ ramified over a reduced plane curve of degree $16$ with a quadruple point and five $(5,5)$-points. The construction of such a curve seems to be a nontrivial computational problem.

We can be more precise about the components of the branch locus $B'+\sum_1^5 A_i\subset W.$

\begin{theorem}\label{components}
Let $S$ be a smooth minimal surface of general type with $p_g=0$ and $K^2=3$ having an involution $i$ such that the bicanonical map of $S$ is not composed with $i.$

With the previous notation, one of the following holds $($here $\Gamma_{a,b}$ denotes a smooth irreducible curve with self-intersection $a$ and genus $b)$: 
\begin{description}
\item[a)] $B'=\Gamma_{-6,0},$ $K_W^2=-1,$ or
\item[b)] $B'=\Gamma_{-6,0}+\Gamma_{-4,0},$ $K_W^2=-2,$ or
\item[c)] $B'=\Gamma_{-2,1}+\Gamma_{-4,0},$ $K_W^2=-1,$ or
\item[d)] $B'=\Gamma_{-2,1}+\Gamma_{-4,0}+\Gamma'_{-4,0},$ $K_W^2=-2,$ or
\item[e)] $B'=\Gamma_{-2,1},$ $K_W^2=0$ $($and ${\rm Kod}(W)=1 )$. 
\end{description}
Moreover, if ${\rm Kod}(W)=1,$ the possibilities for the multiple fibres $m_iF_i$ of the elliptic fibration of $W$ are:\\
$(m_1,m_2,m_3)=(2,2,2)$ or $(m_1,m_2)=(2,3),$ $(2,4)$ or $(3,3).$\\
The corresponding fibrations in $S$ are of genus $3, 7, 5$ or $4,$ respectively.

There is an example for {\rm a)} and {\rm c)}.
\end{theorem}
\begin{remark}
It is immediate from this theorem that the surface $S$ contains at least a smooth rational curve or a smooth elliptic curve. In cases {\rm b)}, {\rm c)} and {\rm d)}, $K_S$ is not ample.
\end{remark}
{\bf Proof of Theorem \ref{components}:}\\
We have $(2K_W+B')B'=4K_WL+4L^2+10=4L(K_W+L)+10=2.$ Since $B\equiv 0\ ({\rm mod}\ 2)$ and $2K_W+B'$ is nef (because $2K_S$ is nef), $B'$ contains an irreducible component $\Gamma$ such that $(2K_W+B')\Gamma=2$ and possibly some components $\Gamma_1,\ldots,\Gamma_l$ such that $(2K_W+B')\Gamma_i=0,$ $i=1,\ldots,l.$
These components are $(-4)$-curves, because $K_V\pi^*(\Gamma_i)=0$ implies that the support of $\pi^*(\Gamma_i)$ is a $(-2)$-curve.

Denote by $\Gamma_V$ the support of $\pi^*(\Gamma).$ One has $K_V\Gamma_V=1$ and then $2g(\Gamma_V)=3+\Gamma_V^2\geq 0.$
The fact $(2K_W+B')(2K_W+B'-3\Gamma)=0$ implies, from the Algebraic Index Theorem, that $(2K_W+B'-3\Gamma)^2\leq 0.$ This gives $\Gamma^2\leq 0,$ thus $\Gamma_V^2=-1$ or $-3$ (equivalently $\Gamma^2=-2$ or $-6$).

We have seen above that $K_W^2\geq -2$ (Proposition \ref{NumR}, d)) and that, if $W$ is birational to an Enriques surface, $K_W^2\leq -1$ ($\overline B$ is singular). Now we claim that $K_W^2\leq -1$ if $W$ is rational. In fact, $-K_W(2K_W+B')=-2$ and $2K_W+B'$ is nef, hence $h^0(W,\mathcal O_W(-K_W))=0$ and then $K_W^2\leq -1$ from the Riemann-Roch Theorem.

Now from $2-2K_W^2=K_WB'=K_W\Gamma+2l=2g(\Gamma)-2-\Gamma^2+2l$ one gets
$$-2K_W^2+4+\Gamma^2=2g(\Gamma)+2l.$$ 
The possibilities allowed by this equation are:
\begin{description}
\item[$\cdot$] $g(\Gamma)=0$, $\Gamma^2=-6$ and $(K_W^2,l)=(-1,0)$ or $(-2,1)$;
\item[$\cdot$] $g(\Gamma)=1,$ $\Gamma^2=-2$ and $(K_W^2,l)=(0,0),$ $(-1,1)$ or $(-2,2)$.
\end{description}

Finally we prove the assertion about the multiple fibres of the elliptic fibration of $W,$ in the case ${\rm Kod}(W)=1.$
The canonical bundle formula (see {\it e.g.} \cite[V. 12.3]{BHPV}) gives $K_P\equiv -F+\sum_1^n(m_i-1)F_i,$ where $m_iF_i\equiv F$ is a multiple fibre of the elliptic fibration of the minimal model $P$ of $W,$ $i=1,\ldots,n.$ Since $K_P\overline B=2,$ we have then
$$\overline BF\left( -1+\sum_1^n\frac{m_i-1}{m_i} \right)=2\ \ {\rm and}\ \ \overline BF\geq 2m_i,\ i=1,\ldots,n.$$
This immediately yields $n\leq 3$ and $n=3\Rightarrow m_1=m_2=m_3=2.$
It is not difficult to see that if $n=2,$ then $(m_1,m_2)=(2,3),$ $(2,4)$ or $(3,3).$

The example for cases a) and c) is given in Section \ref{example}. \Qed

\section{Example}\label{example}

\subsection{Bidouble covers}
A bidouble cover is a finite flat Galois morphism with Galois group $\mathbb Z_2^2.$
Following \cite{Ca} or \cite{Pa}, to define a bidouble cover $V\rightarrow X,$
with $V,$ $X$ smooth surfaces, it suffices to present:
\begin{description}
\item[$\cdot$] smooth divisors $D_1, D_2, D_3\subset X$ with pairwise
transverse intersections and no common intersection;
\item[$\cdot$] line bundles $L_1, L_2, L_3$ such that $2L_g\equiv D_j+D_k$ for
each permutation $(g,j,k)$ of $(1,2,3).$
\end{description}
If ${\rm Pic}(X)$ has no 2-torsion, the $L_i$'s are uniquely determined
by the $D_i$'s.

Let $N:=2K_X+\sum_1^3L_i.$ One has $2K_V\equiv\psi^*\left(N\right)$
and
$$H^0(V,\mathcal O_V(2K_V))\simeq H^0(X,\mathcal O_X(N))\oplus
\bigoplus_{i=1}^3 H^0(X,\mathcal O_X(N-L_i)).$$
The bicanonical map of $V$ is composed with the involution $i_g,$ associated to $L_g,$
if and only if $$h^0(X,\mathcal O_X(2K_X+L_g+L_j))=h^0(X,\mathcal O_X(2K_X+L_g+L_k))=0.$$

For more information on bidouble covers see \cite{Ca} or \cite{Pa}.

\subsection{The construction}\label{construction}

In this section we construct a bidouble cover $V\rightarrow X,$ determined by divisors $D_1,$ $D_2,$ $D_3,$ such that $X$ is a rational surface and the minimal model $S$ of $V$ is a surface of general type with $p_g=0$ and $K^2=3.$
Let $i_g$ be the involution of $V$ corresponding to $2L_g\equiv D_j+D_k$ for each permutation $(g,j,k)$ of $(1,2,3).$ We verify below that the quotients $W_g:=V/{i_g}$ satisfy:
\begin{description}
\item[$\cdot$] $W_1$ is birational to an Enriques surface;
\item[$\cdot$] ${\rm Kod}(W_2)=1,$ $p_g(W_2)=q(W_2)=0;$
\item[$\cdot$] $W_3$ is rational.
\end{description}
Moreover, the surface $S$ has an hyperelliptic fibration of genus $3$ and the bicanonical map $\phi_2$ of $S$ is of degree $2$ onto a rational surface. The map $\phi_2$ is composed with the involution induced by $i_3$ and is not composed with the involutions induced by $i_1$ and $i_2$.
\\\\
{\bf Step 1:} {\em Construction of $S.$}\\
Let $p_0, p_1, p_2\in\mathbb P^2$ be distinct points and $T_1, T_2$ be the lines $p_0p_1, p_0p_2,$ respectively.
In the Appendix we use the Magma functions $LinSys$ and $ParSch$ given in \cite{Ri2} to compute plane curves $C_6$ of degree $6$ and $C_5$ of degree $5$ such that:
\begin{description}
\item[$\cdot$] the singularities of $C_6$ are a double point at $p_0,$ $(2,2)$-points at $p_1, p_2$ tangent to $T_1, T_2$ and a triple point $p_3$ which resolves to a $(2,2)$-point after one blow-up;
\item[$\cdot$] the singularities of $C_5$ are a $(2,2)$-point at $p_1$ tangent to $T_1$ and a $(2,2,2,2)$-point at $p_2$ tangent to $T_2$ such that the intersection number of $C_5$ and $C_6$ at $p_2$ is $12;$
\item[$\cdot$] $C_5$ contains $p_0$ and intersects $C_6$ with multiplicity $7$ at $p_3.$
\end{description}

\begin{figure}[h]

$$
\begin{array}{lll}

\begin{picture}(45,30)
\multiput(0,0)(0,4.25){10}{\line(0,1){2}}
\multiput(-5,20)(4.25,0){12}{\line(1,0){2}}
\qbezier[35](10,40)(10,20)(10,0)
\qbezier[35](25,40)(25,20)(25,0)
\qbezier(18,7)(32,20)(18,33)
\qbezier(32,7)(18,20)(32,33)
\put(-6,33){\makebox(0,0){\tiny $-2$}}
\put(40,16){\makebox(0,0){\tiny $-1$}}
\end{picture}

\hspace{2.5cm}

&

\begin{picture}(45,30)
\multiput(20,0)(0,4.25){10}{\line(0,1){2}}
\qbezier(40,20)(20,20)(0,20)
\qbezier[35](0,5)(20,5)(40,5)
\qbezier[28](7,13)(20,27)(33,13)
\qbezier[28](7,27)(20,13)(33,27)
\put(14,33){\makebox(0,0){\tiny $-1$}}
\end{picture}

\end{array}
$$
\vspace{-15pt}
\caption*{Blow-ups at $p_2$ and $p_3$. \ \ \ 
\begin{picture}(25,5)
\qbezier(0,2.5)(10,2.5)(20,2.5)
\end{picture}
$\widehat{C_5}$. \ \ \ 
\begin{picture}(25,5)
\qbezier[15](0,2.5)(10,2.5)(20,2.5)
\end{picture}
$\widehat{C_6}$.
}
\end{figure}
Let $\mu:X\rightarrow\mathbb P^2$ be the map which resolves the singularities of $C_5+C_6$ and let $E_i,E_i',\ldots$ be the exceptional divisors (with self-intersection $(-1)$) corresponding to the blow-ups at $p_i,$ $i=0,\ldots,3.$ 
Let $T$ denote a general line in $\mathbb P^2$ and let the notation $\widetilde{\cdot}$ denote the total transform $\mu^*(\cdot)$ of a curve.

Let $V\rightarrow X$ be the bidouble cover determined by the divisors
$$
\begin{array}{l}
D_1:=\widetilde{C_5}-E_0-(2E_1+2E_1')-(2E_2+2E_2'+2E_2''+2E_2''')-(E_3+E_3'+E_3''),\\
D_2:=\widetilde{T_1}-E_0-2E_1'+E_3-E_3'+E_3'',\\
D_3:=\widetilde{C_6}+\widetilde{T_2}\\
\ \ \ \ \ \ \ \ -3E_0-(2E_1+2E_1')-(2E_2+4E_2'+0E_2''+2E_2''')-(3E_3+E_3'+3E_3'')
\end{array}
$$
and let $S$ be the minimal model of $V.$
\\\\
{\bf Step 2:} {\em Invariants of $S.$}\\
We have
$$
\begin{array}{l}
L_1\equiv 4\widetilde{T}-2E_0-(E_1+2E_1')-(E_2+2E_2'+E_2''')-(E_3+E_3'+E_3''),\\
L_2\equiv 6\widetilde{T}-2E_0-(2E_1+2E_1')-(2E_2+3E_2'+E_2''+2E_2''')-(2E_3+E_3'+2E_3''),\\
L_3\equiv 3\widetilde{T}-E_0-(E_1+2E_1')-(E_2+E_2'+E_2''+E_2''')-E_3'
\end{array}
$$
and
$$
\begin{array}{l}
K_X+L_1\equiv \widetilde{T}-E_0-E_1'-E_2'+E_2'',\\
K_X+L_2\equiv 3\widetilde{T}-E_0-(E_1+E_1')-(E_2+2E_2'+E_2''')-(E_3+E_3''),\\
K_X+L_3\equiv -E_1'+E_3+E_3''.
\end{array}
$$
One has
$$\chi(\mathcal O_S)=4\chi(\mathcal O_X)+\frac{1}{2}\sum_1^3 L_i(K_X+L_i)=4-1-1-1=1$$
and
$$p_g(S)=p_g(X)+\sum_1^3 h^0(X,\mathcal O_X(K_X+L_i))=0$$
(see the Appendix for the computation of $h^0(X,\mathcal O_X(K_X+L_2))$).
\\\\
{\bf Step 3:} {\em Calculation of $K_S^2.$}\\
Let $N=2K_X+\sum_1^3L_i$.
From the computations in the Appendix we get
$$h^0(V,\mathcal O_V(2K_V))=h^0(X,\mathcal O_X(N))+\sum_{i=1}^3 h^0(X,\mathcal O_X(N-L_i))=4.$$
The surface $V$ contains at least eight $(-1)$-curves (in $\widetilde{T_1}+\widetilde{T_2}$), hence $K_S^2\geq K_V^2+8=N^2+8=1$ and then $S$ is of general type.
Since $h^0(V,\mathcal O_V(2K_V))=h^0(S,\mathcal O_S(2K_S))=K_S^2+1$
(see {\it e.g.} \cite[VII. 5.]{BHPV}), then $K_S^2=3.$
\\\\
{\bf Step 4:} {\em The surface $W_1.$}\\
Let $W_1$ be the double cover of $X$ with branch locus $D_2+D_3.$ It is well known that the smooth minimal model of $W_1$ is an Enriques surface (see {\it e.g.} \cite{CD}).
\\\\
{\bf Step 5:} {\em The surface $W_2.$}\\
Let $W_2$ be the double cover of $X$ with branch locus $D_1+D_3.$
One has $$\chi(\mathcal O_{W_2})=2\chi(\mathcal O_X)+\frac{1}{2}L_2(K_X+L_2)=2-1=1$$ and $$p_g(W_2)=p_g(X)+h^0(X,\mathcal O_X(K_X+L_2))=0.$$
We show in the Appendix that
$$h^0(X,\mathcal O_X(2K_X+2L_2))=1 \ \ \ {\rm and} \ \ \ h^0(X,\mathcal O_X(6K_X+6L_2))=2.$$
This implies ${\rm Kod} (W_2)>0$ and, since
$$h^0(W_2,\mathcal O_{W_2}(2K_{W_2}))=h^0(X,\mathcal O_X(2K_X+L_2))+h^0(X,\mathcal O_X(2K_X+2L_2))=1,$$
$W_2$ is not of general type (see {\it e.g.} \cite[VII. 5.]{BHPV}).
This way ${\rm Kod} (W_2)=1.$
\\\\
{\bf Step 6:} {\em The surface $W_3.$}\\
Let $\rho:W_3\rightarrow X$ be the double cover with branch locus $D_1+D_2.$
The pencil of conics tangent to the lines $T_1, T_2$ at $p_1, p_2$ lifts to a rational fibration of $W_3$ (and lifts to a genus $3$ fibration of $S$). Since
$$\chi(\mathcal O_{W_3})=2\chi(\mathcal O_X)+\frac{1}{2}L_3(K_X+L_3)=2-1=1$$
and $$p_g(W_3)=p_g(X)+h^0(X,\mathcal O_X(K_X+L_3))=0,$$
then $W_3$ is a rational surface.
\\\\
{\bf Step 7:} {\em Bicanonical map.}\\
As computed in the Appendix, one has
$$
\begin{array}{l}
h^0(X,\mathcal O_X(2K_X+L_1+L_2))=1,\\
h^0(X,\mathcal O_X(2K_X+L_1+L_3))=0,\\
h^0(X,\mathcal O_X(2K_X+L_2+L_3))=0,
\end{array}
$$
hence the bicanonical map $\phi_2'$ of $V$ is not composed with the involutions $i_1$ and $i_2$ and is composed with $i_3.$
Let $\psi_1:V\rightarrow W_3$ be the double cover corresponding to $i_3$ and let
$\psi_2:W_3\rightarrow\mathbb P^3$ be the map induced by
$$H^0(X,\mathcal O_X(\rho^*(2K_X+L_1+L_2+L_3)))\oplus H^0(X,\mathcal O_X(\rho^*(2K_X+L_1+L_2)+R)),$$
where $R$ is the ramification divisor of the map $\rho:W_3\rightarrow X$ defined above. We have $$\phi_2'=\psi_1\circ\psi_2.$$
It is shown in the Appendix that the degree of $\psi_2(W_3)$ is $6.$ Since $(2K_S)^2=12,$ this implies that the bicanonical map of $S$ is of degree $2.$

\appendix
\section*{Appendix: Magma code}

Here the {\em Computational Algebra System Magma} (\cite{BCP}) is used to perform some calculations.\\
We use the following Magma functions, given in \cite{Ri2}: $LinSys,$ which computes linear systems of plane curves with non-ordinary singularities, and $ParSch,$ whose output is a scheme which parametrizes given degree plane curves with given singularities.\\

\noindent
{\bf 1)} First we compute the curves $C_5$ and $C_6$ referred in Section \ref{construction}.

\small
\begin{verbatim}
K:=Rationals();
A<x,y>:=AffineSpace(K,2);
L:=[LinearSystem(A,6),LinearSystem(A,5),LinearSystem(A,3),\
	LinearSystem(A,2),LinearSystem(A,2),LinearSystem(A,1),LinearSystem(A,1)];
P:=[A![0,0],A![0,1],A![1,0],A![1,1]];
M:=[[[2],[2,2],[2,2,1,1],[3,2,2]],\
    [[1],[2,2],[2,2,2,2],[1,1,1]],\
    [[1],[1,1],[1,1,1,1],[1,1,1]],\
    [[0],[1,1],[1,1,0,0],[1,1,0]],\
    [[0],[1,1],[1,1,1,1],[0,0,0]],\
    [[0],[0,0],[1,1,1,0],[0,0,0]],\
    [[0],[0,0],[0,0,0,0],[1,1,1]]
];
T:=[[],[[0,1]],[[1,0],[],[]],[[],[]]];
\end{verbatim}
\normalsize
We want to compute points (some infinitely near) such that:\\
$\cdot$ the sets of elements of $L[1], L[2]$ which have singularities, at those points, of multiplicities given by $M[1], M[2]$ are non-empty;\\
$\cdot$ the five sets of elements of $L[3],\ldots,L[7]$ of curves with singularities, at those points, of multiplicities given by $M[3],\ldots,M[7]$ are empty.\\
This last step is needed in order to obtain a non-reduced curve.
The following gives a scheme which parametrizes such curves.
\small
\begin{verbatim}
S:=ParSch(L,P,M,T,[],[],5);
\end{verbatim}
\normalsize
This scheme is zero dimensional. We compute a point in $S$
\small
\begin{verbatim}
PointsOverSplittingField(S);
\end{verbatim}
\normalsize
and we use the function $LinSys$ to compute the reduced curves $C_5$ and $C_6.$
\small
\begin{verbatim}
R<r1>:=PolynomialRing(Rationals());
K<r1>:=NumberField(r1^2 + 1496/675*r1 + 10976/625);
A<x,y>:=AffineSpace(K,2);
L5:=LinearSystem(A,5);L6:=LinearSystem(A,6);
P:=[A![0,0],A![0,1],A![1,0],A![1,1]];
M5:=[[1],[2,2],[2,2,2,2],[1,1,1]];M6:=[[2],[2,2],[2,2,1,1],[3,2,2]];
T:=[[],[[0,1]],[[1,0],[1,15/61*r1 + 443/2745],[1,35/366*r1 - 10802/24705]],\
[[1,-405/1708*r1 - 1171/2135],[1,-21465/95648*r1 - 23559/59780]]];
J5:=LinSys(L5,P,M5,T);C_5:=Curve(A,Sections(J5)[1]);
J6:=LinSys(L6,P,M6,T);C_6:=Curve(A,Sections(J6)[1]);
\end{verbatim}
\normalsize
The equations of $C_5$ and $C_6$ are, in affine space:
\footnotesize
\begin{verbatim}
3660*x^5+(-900*r+1341)*x^4*y-14640*x^4+(-3550*r-12858)*x^3*y^2+(4500*r-300)*x^3*y
+21960*x^3+(-1550*r-14313)*x^2*y^3+(7800*r+34128)*x^2*y^2+(-6300*r-4338)*x^2*y
-14640*x^2+(1350*r-8874)*x*y^4+29280*x*y^3+(-4050*r-28278)*x*y^2+(2700*r+4212)*x*y
+3660*x-915*y^5+3660*y^4-5490*y^3+3660*y^2-915*y
\end{verbatim}
\normalsize
and
\footnotesize
\begin{verbatim}
35882945*x^6+(-161700*r+36034208)*x^5*y-143531780*x^5+(-12929700*r-26583872)*x^4
*y^2+(13414800*r-81518752)*x^4*y+215297670*x^4+(4648050*r-9108022)*x^3*y^3
+(16563300*r+71383788)*x^3*y^2+(-21696450*r+45826858)*x^3*y-143531780*x^3
+(12738300*r-1064672)*x^2*y^4+(-34772700*r+20345388)*x^2*y^3+(18400800*r-64080632)
*x^2*y^2+(3795300*r+8765708)*x^2*y+35882945*x^2+(666300*r+36857408)*x*y^5
+(-14737200*r-109507552)*x*y^4+(32123550*r+99334858)*x*y^3+(-22700700*r-17576692)
*x*y^2+(4648050*r-9108022)*x*y+(166050*r+477113)*y^6+(-664200*r-1908452)*y^5
+(996300*r+2862678)*y^4+(-664200*r-1908452)*y^3+(166050*r+477113)*y^2
\end{verbatim}
\normalsize
with
\footnotesize
\verb-r^2+1496/675*r+10976/625=0-.
\normalsize
\\\\
\noindent
{\bf 2)} From Section \ref{construction} one has:
\small
$$
\begin{array}{l}
2K_X+\sum_1^3 L_i\equiv 7\widetilde{T}-3E_0-(2E_1+4E_1')-(2E_2+4E_2'+2E_2''')-(E_3+E_3'+E_3''),\\
K_X+L_2\equiv 3\widetilde{T}-E_0-(E_1+E_1')-(E_2+2E_2'+E_2''')-(E_3+E_3''),\\
2K_X+2L_2\equiv 6\widetilde{T}-2E_0-(2E_1+2E_1')-(2E_2+4E_2'+2E_2''')-(2E_3+2E_3''),\\
6K_X+6L_2\equiv 18\widetilde{T}-6E_0-(6E_1+6E_1')-(6E_2+12E_2'+6E_2''')-(6E_3+6E_3''),\\
2K_X+L_1+L_2\equiv 4\widetilde{T}-2E_0-(E_1+2E_1')-(E_2+3E_2'-E_2''+E_2''')-(E_3+E_3''),\\
2K_X+L_1+L_3\equiv \widetilde{T}-E_0-2E_1'-E_2'+E_2''+E_3+E_3'',\\
2K_X+L_2+L_3\equiv 3\widetilde{T}-E_0-(E_1+2E_1')-(E_2+2E_2'+E_2''').
\end{array}
$$
Below we compute the dimension of the first cohomology group $h^0(X,\mathcal O_X(\cdot))$ for each of these divisors
(it is immediate that $h^0(X,\mathcal O_X(2K_X+L_1+L_3))=0$). We obtain $3, 0, 1, 2, 1, 0, 0,$ respectively.
\small
\begin{verbatim}
M:=[[[3],[2,4],[2, 4,0,2],[1,1,1]],
    [[1],[1,1],[1, 2,0,1],[1,0,1]],
    [[2],[2,2],[2, 4,0,2],[2,0,2]],
    [[6],[6,6],[6,12,0,6],[6,0,6]],
    [[2],[1,2],[1, 3,0,0],[1,0,1]],
    [[1],[1,2],[1, 2,0,1],[0,0,0]]];
d:=[7,3,6,18,4,3];
J:=[LinSys(LinearSystem(A,d[i]),P,M[i],T):i in [1..6]];
[#Sections(J[i]):i in [1..6]];
\end{verbatim}
\normalsize
\noindent
{\bf 3)} Now we describe how to compute the degree of the scheme $\psi_2(W_3)$ referred in Section \ref{construction}, Step $7$. The complete code is available at\\
\verb+http://home.utad.pt/~crito/magma_code.html+
\\\\
Let $f_6$ be the defining equation of the curve $D_1+D_2,$ $f_4$ be the equation of the unique effective plane curve corresponding to $2K_X+L_1+L_2$ and let $J_7$ be the linear system of plane curves corresponding to $2K_X+L_1+L_2+L_3.$ We define (a singular model of) $W_3$ in a weighted projective space and we define the map $\psi_2:W_3\rightarrow\mathbb P^3.$
\small
\begin{verbatim}
WP<w,x,y,z>:=ProjectiveSpace(K,[3,1,1,1]);
W3:=Scheme(WP,w^2-f6);
P3:=ProjectiveSpace(K,3);
psi2:=map<W3->P3|(Sections(J7) div (x*y)) cat [w*(f4 div (x*y))]>;
\end{verbatim}
\normalsize
We want to compute
\small
\begin{verbatim}
Degree(psi2(W3));
\end{verbatim}
\normalsize
but $6$ GB of computer memory are not enough for this task.
Thus we compute the degree of the intersection of two hyperplane sections of $\psi_2(W_3)$.
\small
\begin{verbatim}
Degree(psi2(Scheme(W3,[Sections(J7)[1],Sections(J7)[2]])));
\end{verbatim}
\normalsize
We obtain degree $6$.

\bibliography{ReferencesRito}

\begin{thebibliography}{BHPV}
\expandafter\ifx\csname urlstyle\endcsname\relax
  \expandafter\ifx\csname doi\endcsname\relax
  \def\doi#1{doi:\discretionary{}{}{}#1}\fi \else
  \expandafter\ifx\csname doi\endcsname\relax
  \def\doi{doi:\discretionary{}{}{}\begingroup \urlstyle{rm}\Url}\fi \fi

\bibitem[BHPV]{BHPV}
W.~Barth, K.~Hulek, C.~Peters and A.~Van~de Ven, {\em {Compact complex
  surfaces. $2$nd enlarged ed.}\/}, {Ergebnisse der Mathematik und ihrer
  Grenzgebiete. 3. Folge 4. Berlin: Springer. xii, 436~p.} (2004).

\bibitem[BP]{BP}
I.~Bauer and R.~Pignatelli, {\em The classification of minimal product-quotient
  surfaces with $p_g=0$\/} (2010), arXiv:1006.3209v1 [math.AG].

\bibitem[Be]{Be}
A.~Beauville, {\em Surfaces alg\' ebriques complexes\/}, vol.~54, Ast\' erisque
  (1978).

\bibitem[BCP]{BCP}
W.~Bosma, J.~Cannon and C.~Playoust, {\em The Magma algebra system. I. The user
  language.\/}, J. Symbolic Comput., {\bf 24} (1997), no. 3--4, 235--265.

\bibitem[Bu]{Bu}
P.~Burniat, {\em {Sur les surfaces de genre $P_{12}>0$}\/}, Ann. Mat. Pura
  Appl., IV. Ser., {\bf 71} (1966), 1--24.

\bibitem[CCM]{CCM}
A.~Calabri, C.~Ciliberto and M.~{Mendes Lopes}, {\em Numerical Godeaux surfaces
  with an involution\/}, Trans. Am. Math. Soc., {\bf 359} (2007), no.~4,
  1605--1632.

\bibitem[Ca]{Ca}
F.~Catanese, {\em {Singular bidouble covers and the construction of interesting
  algebraic surfaces}\/}, {Contemp. Math. 241, Am. Math. Soc, 97--120} (1999).

\bibitem[CFM]{CFM}
C.~Ciliberto, P.~Francia and M.~{Mendes Lopes}, {\em Remarks on the bicanonical
  map for surfaces of general type\/}, Math. Z., {\bf 224} (1997), no.~1,
  137--166.

\bibitem[CM]{CM}
C.~Ciliberto and M.~{Mendes Lopes}, {\em On surfaces with $p_g=q=2$ and
  non-birational bicanonical map\/}, Adv. Geom., {\bf 2} (2002), no.~3,
  281--300.

\bibitem[CD]{CD}
F.~Cossec and I.~Dolgachev, {\em {Enriques surfaces. I.}\/}, {Progress in
  Mathematics, 76. Boston, MA etc.: Birkh\"auser Verlag. ix, 397 p. DM 96.00 }
  (1989).

\bibitem[In]{In}
M.~Inoue, {\em {Some new surfaces of general type}\/}, Tokyo J. Math., {\bf 17}
  (1994), no.~2, 295--319.

\bibitem[Ke]{Ke}
J.~H. Keum, {\em Some new surfaces of general type with $p_g=0$\/}, preprint,
  1988.

\bibitem[MP1]{MP1}
M.~Mendes~Lopes and R.~Pardini, {\em {The bicanonical map of surfaces with
  $p\sb g=0$ and $K\sp 2\geq 7$.}\/}, Bull. Lond. Math. Soc., {\bf 33} (2001),
  no.~3, 265--274.

\bibitem[MP2]{MP2}
M.~Mendes~Lopes and R.~Pardini, {\em {A survey on the bicanonical map of
  surfaces with $p_g=0$ and $K^2\geq 2$.}\/}, {Beltrametti, Mauro C. (ed.) et
  al., Algebraic geometry. A volume in memory of Paolo Francia. Berlin: de
  Gruyter. 277-287 (2002).} (2002).

\bibitem[MP3]{MP3}
M.~Mendes~Lopes and R.~Pardini, {\em {Enriques surfaces with eight nodes.}\/},
  Math. Z., {\bf 241} (2002), no.~4, 673--683.

\bibitem[MP4]{MP4}
M.~Mendes~Lopes and R.~Pardini, {\em {The bicanonical map of surfaces with
  $p_g=0$ and $K^2\geq 7$. II.}\/}, Bull. Lond. Math. Soc., {\bf 35} (2003),
  no.~3, 337--343.

\bibitem[MP5]{MP5}
M.~Mendes~Lopes and R.~Pardini, {\em {A new family of surfaces with $p_{g}=0$
  and $K^2=3$.}\/}, Ann. Sci. Éc. Norm. Supér. (4), {\bf 37} (2004), no.~4,
  507--531.

\bibitem[MP6]{MP6}
M.~Mendes~Lopes and R.~Pardini, {\em {Surfaces of general type with $p_g=0$,
  $K^2=6$ and non birational bicanonical map.}\/}, Math. Ann., {\bf 329}
  (2004), no.~3, 535--552.

\bibitem[Na]{Na}
D.~Naie, {\em {Surfaces d'Enriques et une construction de surfaces de type
  g\'en\'eral avec $p\sb g=0$}\/}, Math. Z., {\bf 215} (1994), no.~2, 269--280.

\bibitem[Pa]{Pa}
R.~Pardini, {\em Abelian covers of algebraic varieties\/}, J. Reine Angew.
  Math., {\bf 417} (1991), 191--213.

\bibitem[PPS2]{PPS2}
H.~Park, J.~Park and D.~Shin, {\em A complex surface of general type with
  $p_{g}=0,$ $K^{2}=3$ and $H_1=\mathbb Z/2\mathbb Z$\/} (2008),
  arXiv:0803.1322v2 [math.AG].

\bibitem[PPS1]{PPS1}
H.~Park, J.~Park and D.~Shin, {\em {A simply connected surface of general type
  with $p_{g}=0$ and $K^{2}=3$.}\/}, Geom. Topol., {\bf 13} (2009), no.~2,
  743--767.

\bibitem[Ri1]{Ri1}
C.~Rito, {\em {Involutions on surfaces with $p_g = q = 1$}\/}, Collect. Math.,
  {\bf 61} (2009), no.~1, 81--106.

\bibitem[Ri2]{Ri2}
C.~Rito, {\em On the computation of singular plane curves and quartic
  surfaces\/} (2010), arXiv:0906.3480v3 [math.AG].

\end{thebibliography}

\bigskip
\bigskip

\noindent Carlos Rito
\\ Departamento de Matem\' atica
\\ Universidade de Tr\' as-os-Montes e Alto Douro
\\ 5001-801 Vila Real
\\ Portugal
\\\\
\noindent {\it e-mail:} crito@utad.pt

\end{document}